\documentclass[11pt,twoside]{article}
\usepackage{fancyhdr}
\usepackage[colorlinks,citecolor=blue,urlcolor=blue,linkcolor=blue,bookmarks=false]{hyperref}
\usepackage{amsfonts,epsfig,graphicx}
\usepackage{afterpage}
\usepackage{amsmath,amssymb,amsthm} 
\usepackage{fullpage}
\usepackage[T1]{fontenc} 
\usepackage{epsf} 
\usepackage{graphics} 
\usepackage{amsfonts,amsmath}
\usepackage[sort,numbers]{natbib} 
\usepackage{psfrag,xspace}
\usepackage{color,etoolbox}

\usepackage{amsmath,amssymb,fullpage,graphicx}

\let\widehat\widehat

\newtheorem{theorem}{Theorem}
\newtheorem{lemma}[theorem]{Lemma}

\newtheorem{proposition}[theorem]{Proposition}
\newtheorem{remark}[theorem]{Remark}
\DeclareMathOperator*{\argmax}{argmax}

\begin{document}

\begin{center} {\LARGE{\bf{Universal Inference}}}
\\

\vspace*{.3in}

{\large{
\begin{tabular}{ccccc}
Larry Wasserman & Aaditya Ramdas & Sivaraman Balakrishnan \\
\end{tabular}

\vspace*{.1in}

\begin{tabular}{ccc}
Department of Statistics and Data Science \\
Machine Learning Department \\
\end{tabular}

\begin{tabular}{c}
Carnegie Mellon University, \\
Pittsburgh, PA 15213.
\end{tabular}

\vspace*{.2in}

\begin{tabular}{c}
{\texttt{\{larry, aramdas, siva\}@stat.cmu.edu}}
\end{tabular}
}}

\vspace*{.2in}

\today
\vspace*{.2in}
\begin{abstract}
We propose a general method for
constructing
confidence sets and hypothesis tests that have finite-sample guarantees
without regularity conditions. We refer to such procedures as ``universal.''
The method is very simple and is based
on a modified version of the usual
likelihood ratio statistic,
that we call ``the split likelihood ratio test'' (split LRT) statistic.
The (limiting) null distribution of the classical likelihood ratio statistic
is often intractable when used to test composite null hypotheses in irregular statistical
models.
Our method is especially appealing for statistical inference
in these complex setups.
The method we suggest works for any parametric model
and also for some nonparametric models, 
as long as computing 
a maximum likelihood estimator (MLE) is feasible under the null. 
Canonical examples arise in mixture modeling and
shape-constrained inference, for which
constructing tests and confidence sets has been notoriously difficult.
We also develop various extensions of our basic methods. 
We show that in settings when computing the MLE is hard, for the purpose
of constructing valid tests and intervals, it is sufficient to upper bound the maximum likelihood.
We investigate some conditions under which our methods yield valid inferences under model-misspecification.
Further, the split LRT can be used with profile likelihoods to deal with nuisance parameters, 
and it can also be run sequentially to yield anytime-valid $p$-values and confidence sequences. Finally, when combined 
with the method of sieves, it can be used to perform model selection with nested model classes. 
\end{abstract}
\end{center}

\section{Introduction}
The foundations of statistics are built on 
a variety of generally applicable principles for parametric estimation and inference. 
In parametric statistical
models, the likelihood ratio test, and confidence intervals obtained from asymptotically Gaussian estimators, are the workhorse
inferential tools for constructing hypothesis tests and confidence intervals.
Often, the validity of these methods relies on large sample asymptotic theory
and requires that the statistical model
satisfy certain regularity conditions; see Section~\ref{section::sanity} for precise definitions.
When these conditions do not hold,
there is no general method for statistical inference, and these settings are typically
considered in an ad-hoc manner.
Here, we introduce
a universal method which yields
tests and confidence sets
for any statistical model
and has finite-sample guarantees.

We begin with some terminology.
A parametric statistical model is a collection of distributions
$\{P_\theta:\ \theta\in\Theta\}$
for an arbitrary set $\Theta$.
When the aforementioned regularity conditions hold,
there are many methods for inference.
For example, if $\Theta\subseteq\mathbb{R}^d$, 
the set
\begin{align}\label{eq:LRT-set}
A_n = \left\{\theta:\ 2 \log \frac{{\cal L}(\widehat\theta)}{{\cal L}(\theta)} \leq c_{\alpha,d}\right\}
\end{align}
is the likelihood ratio confidence set,
where $c_{\alpha,d}$ is the upper $\alpha$-quantile of a $\chi^2_d$ distribution, ${\cal L}$ is the likelihood function
and $\widehat\theta$ is the maximum likelihood estimator (MLE).
It satisfies the asymptotic coverage guarantee
\begin{align*}
P_{\theta^*}(\theta^* \in A_n)\to 1-\alpha
\end{align*}
as $n\to\infty$, 
where $P_{\theta^*}$ denotes the unknown true data generating distribution.

Constructing tests and confidence intervals for 
irregular models --- where the regularity conditions do not hold ---
is very difficult \cite{drton2009likelihood}.
An example is
mixture models.
In this case 
we observe $Y_1,\ldots, Y_n \sim P$ and we want to test
\begin{align}
\label{eqn:mixture}
H_0: ~~P\in {\cal M}_{k_0}~~\text{versus}~~H_1: ~~P\in {\cal M}_{k_1},
\end{align}
where ${\cal M}_k$ denotes the set of mixtures of $k$ Gaussians,
with an appropriately restricted parameter space $\Theta$ (see for instance \cite{redner1981}),
and with $k_0 < k_1$.
Finding a test that provably controls the type I error at a given level
has been elusive.
A natural candidate is to base the test on the 
likelihood ratio statistic
but this turns out to have an intractable
limiting distribution
\cite{dacunha1997testing}.
As we discuss further in Section~\ref{sec:examples}, developing practical, simple tests for this pair of hypotheses is an active area of research \cite[and references therein]{chen2009hypothesis,mclachlan1987bootstrapping,purvasha2019gaussian}. 
However, it is possible that we may be able to compute an MLE using variants of the EM algorithm.
In this paper,
we show that there is a remarkably simple test based on the MLE with
guaranteed finite-sample control of the type I error.
Similarly, we construct a confidence set for the parameters of a mixture model
with guaranteed finite-sample coverage.
These tests and confidence sets can in fact be used for any model.
In regular statistical models (those for which the usual LRT is well-behaved), our 
methods may not be optimal, though we do not yet fully understand how close to optimal they are beyond special cases (uniform, Gaussian). Our test is most useful in irregular (or singular) models
for which valid tests are not known, or require many assumptions. 
Going beyond parametric models, we show 
that our methods can be used
for several nonparametric models as well, and has a natural sequential analog.

\section{Universal Inference}

Let
$Y_1,\ldots, Y_{2n}$ be an iid sample from
a distribution $P_{\theta^*}$ which belongs to a collection
$(P_\theta:\ \theta\in\Theta)$.
Note that $\theta^*$ denotes the true value of the parameter.
Assume that each distribution $P_\theta$ has a density $p_\theta$
with respect to some underlying measure $\mu$ (for instance the Lebesgue or counting measure).

\paragraph{A universal confidence set.} We construct a confidence set for $\theta^*$ by
first splitting the data into two groups
$D_0$ and $D_1$.
For simplicity, we take each group to
be of the same size $n$ but this is not necessary.
Let $\widehat\theta_1$ be \emph{any} estimator
constructed from $D_1$; this can be the MLE, a Bayes estimator that utilizes prior knowledge, a robust estimator, etc.
Let
\[
{\cal L}_0(\theta) = \prod_{i\in D_0}p_\theta(Y_i)
\]
denote the likelihood function based on $D_0$. We define the
{\em split likelihood ratio statistic (split LRS)} as
\begin{equation}
T_n(\theta) = \frac{ {\cal L}_0(\widehat\theta_1)}{ {\cal L}_0(\theta)}.
\end{equation}
Then, the universal confidence set is
\begin{equation}
\label{eqn:universal}
C_n = \Biggl\{ \theta \in \Theta :\ T_n(\theta) \leq \frac{1}{\alpha} \Biggr\}.
\end{equation}
Similarly, define the {\em crossfit LRS} as
\begin{equation}
S_n(\theta) = (T_n(\theta) + T_n^{\text{swap}}(\theta))/2,
\end{equation}
where $T_n^{\text{swap}}$ is formed by calculating $T_n$ after swapping the roles of $D_0$ and $D_1$.
We can also define $C_n$ with $S_n$ in place of $T_n$.

\begin{theorem}\label{thm:universal-set}
$C_n$ is a finite-sample valid $(1-\alpha)$ confidence set for $\theta^*$, meaning that $P_{\theta^*}(\theta^* \in C_n) \geq 1-\alpha.$
\end{theorem}

If we did not split the data and $\widehat \theta_1$ was the MLE, then $T_n(\theta)$ would be the usual likelihood ratio statistic and we would typically approximate its distribution using
an asymptotic argument.
For example, as mentioned earlier,
in regular models, -2 times the log likelihood ratio statistic
has, asymptotically, a $\chi^2_d$ distribution.
But, in irregular models this strategy can fail.
Indeed, finding or approximating the distribution of the likelihood ratio statistic
is highly nontrivial in irregular models.
The split LRS avoids these complications.

Now we explain why $C_n$ has coverage at least $1-\alpha$, as claimed by the above theorem. We prove it for the version using $T_n$, but the proof for $S_n$ is identical.
Consider any fixed $\psi\in\Theta$ and let $A$ denote the support of $P_{\theta^*}$.
Then, 
\begin{align*}
\mathbb{E}_{\theta^*} \left[\frac{{\cal L}_0(\psi)}{{\cal L}_0(\theta^*)}\right] &=
\mathbb{E}_{\theta^*} \left[
\frac{ \prod_{i\in D_0}p_\psi(Y_i)}{ \prod_{i\in D_0}p_{\theta^*}(Y_i)}\right]\\
&\hspace{-1.7cm}=
\int_A 
\frac{ \prod_{i\in D_0}p_\psi(y_i)}{ \prod_{i\in D_0}p_{\theta^*}(y_i)}
\prod_{i\in D_0}p_{\theta^*}(y_i) ~ d y_1 \cdots d y_n\\
&\hspace{-1.7cm}=
\int_A  \prod_{i\in D_0}p_\psi(y_i) d y_1 \cdots d y_n \leq 
 \prod_{i\in D_0} \left[ \int p_\psi(y_i) d y_i \right] = 1.
\end{align*}
Since $\widehat\theta_1$ is fixed when we condition on $D_1$,
we have
\begin{equation}\label{eq:unit-expectation}
\mathbb{E}_{\theta^*}[ T_n(\theta^*) \,| \, D_1 ] = \mathbb{E}_{\theta^*}\left[ \frac{{\cal L}_0(\widehat\theta_1)}{{\cal L}_0(\theta^*)} \,\Biggm| \, D_1 \right] \leq 1.
\end{equation}
Now, using Markov's inequality, 
\begin{align}
P_{\theta^*}(\theta^*\notin C_n) &=
P_{\theta^*}\left( T_n(\theta^*) > \frac{1}{\alpha}\right)\leq
  \alpha
\mathbb{E}_{\theta^*} [ T_n(\theta^*)] \label{eq:markov}\\
&\hspace{-1.8cm}=
\alpha
\mathbb{E}_{\theta^*}\left[ \frac{{\cal L}_0(\widehat\theta_1)}{{\cal L}_0(\theta^*)}\right] =
\alpha\mathbb{E}_{\theta^*}\left(\mathbb{E}_{\theta^*}\left[ \frac{{\cal L}_0(\widehat\theta_1)}{{\cal L}_0(\theta^*)} \,\Biggm| \, D_1 \right]\right)
\leq \alpha.
\nonumber
\end{align}

\begin{remark}\label{rem:nonparametric}
The parametric setup adopted above generalizes easily to nonparametric settings as long as we can calculate a likelihood. 
For a collection of densities $\mathcal{P}$, and a true density $p^* \in \mathcal{P}$, suppose 
we use $D_1$ to identify $\widehat p_1 \in \mathcal{P}$, and $D_0$ to calculate 
\begin{align*}
T_n(p) = \prod_{i\in D_0} \frac{ \widehat p_1(Y_i) }{ p(Y_i)}.
\end{align*}
We then define,
$C_n := \{p \in \mathcal{P}: T_n(p) \leq 1/\alpha\},$
and our previous argument ensures that $P_{p^*}(p^* \in C_n) \geq 1 - \alpha$.
\end{remark}

\paragraph{A universal hypothesis test.}
Now we turn to hypothesis testing.
Let $\Theta_0\subset\Theta$ be a possibly composite null set and consider testing
\begin{equation}\label{eq:test}
H_0: \theta^* \in\Theta_0\ \ \ {\rm versus}\ \ \ \theta^* \notin\Theta_0.
\end{equation}
The alternative above can be replaced by $\theta^* \in \Theta_1$ for any $\Theta_1 \subseteq \Theta$ or by $\theta^* \in \Theta_1 \backslash \Theta_0$. 
One way to test this hypothesis is based on the universal confidence set in~\eqref{eqn:universal}. We simply reject the null hypothesis if $C_n \bigcap \Theta_0 = \emptyset.$ 
It is straightforward to see that if this test makes a type I error then the universal confidence set must fail to cover $\theta^*$, and so the type I error of this test
is  at most $\alpha$.

We present an alternative method that is often computationally (and possibly statistically) more attractive. 
Let $\widehat\theta_1$ be any estimator constructed from $D_1$,
and let \[
\widehat\theta_{0}:= \argmax_{\theta \in \Theta_0} {\cal L}_0(\theta)\] 
be the MLE under $H_0$ constructed from $D_0$.
Then the universal test, which we call the \emph{split likelihood ratio test (split LRT)}, is defined as:
\begin{equation}\label{eq:splitLRT}
\text{reject $H_0$ if } U_n > 1/\alpha, \text{ where } U_n = \frac{{\cal L}_0(\widehat\theta_1)}{{\cal L}_0(\widehat\theta_{0})}.
\end{equation}
Similarly, we can define the \emph{crossfit LRT} as
\begin{equation}\label{eq:crossfitLRT}
\text{reject $H_0$ if } W_n > 1/\alpha, \text{ where } W_n = \frac{U_n + U^{\text{swap}}_n}2,
\end{equation}
where as before, $U^{\text{swap}}_n$ is calculated like $U_n$ after swapping the roles of $D_0$ and $D_1$.

\begin{theorem}
The split and crossfit LRTs  control the type~I error at $\alpha$, i.e. $\sup_{\theta^* \in \Theta_0} P_{\theta^*}(U_n > 1/\alpha) \leq \alpha$.
\end{theorem}

The proof is straightforward. We prove it for split LRT, but once again the crossfit proof is identical.
Suppose that $H_0$ is true and
$\theta^*\in\Theta_0$ is the true parameter.
By Markov's inequality, the type I error is
\begin{align*}
P_{\theta^*}( U_n > 1/\alpha) &=
P_{\theta^*} \left(
 {\cal L}_0(\widehat\theta_1)/{\cal L}_0(\widehat\theta_{0}) > 1/ \alpha\right)\\
& \hspace{-1cm} \leq  \alpha
\mathbb{E}_{\theta^*}\left[\frac{{\cal L}_0(\widehat\theta_1)}{{\cal L}_0(\widehat\theta_{0})}\right]
 \stackrel{\text{(i)}}{\leq} \alpha
\mathbb{E}_{\theta^*}\left[\frac{ {\cal L}_0(\widehat\theta_1)}{{\cal L}_0(\theta^*)}\right] \stackrel{\text{(ii)}}{\leq} \alpha.
\end{align*}
Above, inequality~(i) uses the fact that
${\cal L}_0(\widehat\theta_{0}) \geq {\cal L}_0(\theta^*)$
which is true when $\widehat\theta_{0}$ is the MLE, and inequality~(ii) 
follows by conditioning on $D_1$ as argued earlier in \eqref{eq:markov}.

\begin{remark}\label{rem:nonparametric2}
We may drop the use of $\Theta, \Theta_0, \Theta_1$ above, 
and extend the split LRT to a general nonparametric setup.
Both tests can be used to test any null $H_0: p^* \in \mathcal{P}_0$ against any alternative $H_1: p^* \in \mathcal{P}_1$. Importantly, no parametric assumption is needed on $\mathcal{P}_0,\mathcal{P}_1$, and no relationship is imposed whatsoever between $\mathcal{P}_0,\mathcal{P}_1$. As before, use $D_1$ to identify $\widehat p_1 \in \mathcal{P}_1$, use $D_0$ to calculate the MLE $\widehat p_0 \in \mathcal{P}_0$, and define $U_n = \prod_{i\in D_0} \frac{ \widehat p_1(Y_i) }{ \widehat p_0(Y_i)}$.
\end{remark}

We call these procedures {\em universal} to mean that
they are valid in finite-samples with
no regularity conditions.
Constructions like this are reminiscent of
ideas used in sequential settings where an estimator is computed from past data
and the likelihood is evaluated on current data; we expand on this in Section~\ref{sec:seq}.

We note in passing that another universal set is the following.
Define
$
C = \big\{ \theta: \int_\Theta {\cal L}(\psi) d\Pi(\psi)/{\cal L}(\theta) \leq 1/\alpha \big\},
$
where ${\cal L}$ is the full likelihood (from all the data) and $\Pi$ is any prior.
This is also has the same coverage guarantee but
requires specifiying a prior and doing an integral.
In irregular or nonparametric models, the integral will typically
be intractable.

\paragraph{Perspective: Poor man's Chernoff bound.} At first glance, the reader may worry that
 Markov's inequality seems like a weak tool to use, resulting in an underpowered conservative test or
 confidence interval.
 However, this is not the right perspective. One should really view our proof as using a ``poor man's Chernoff bound''.

 For a regular model, we would usually compare the log-likelihood ratio to the $(1-\alpha)$-quantile of a chi-squared distribution (with degrees of freedom related to the difference in dimensionality of the null and alternate models). Instead, we compare the log-split-likelihood ratio to $\log(1/\alpha)$, which
 scales like the $(1-\alpha)$-quantile of a chi-squared distribution with one degree of freedom. 

In any case, instead of finding the asymptotic distribution of $\log
U_n$ (usually having a moment generating function,
like a chi-squared), our proof should be interpreted as using the
simpler but nontrivial fact that $\mathbb{E}_{\theta^*}[e^{\log(U_n)}]
\leq 1$. Hence we are really using the fact that $\log U_n$ has an
exponential tail, just as an asymptotic argument would. 

A true Chernoff-style bound for a chi-squared random
variable would have bounded $\mathbb{E}_{\theta^*}[e^{a\log(U_n)}]$ by
an appropriate function of $a$, and then optimized over the choice of $a > 0$ 
to obtain a tight bound. Our methods correspond to choosing $a = 1$, leading 
us to call the technique a poor man's Chernoff bound.  The key point is that our methods
should be viewed as using Markov's inequality on the exponential of the random 
variable of interest. 



\paragraph{Perspective: In-sample versus out-of-sample likelihood.}
We may rewrite the universal set as
\[
C_n = \left\{ \theta \in \Theta :\ 2 \log \frac{ {\cal L}_0(\widehat\theta_1)}{ {\cal L}_0(\theta)} \leq 2 \log (1/\alpha) \right\}.
\]
For a regular model, it is natural to compare the above expression to
the usual LRT-based set $A_n$ from \eqref{eq:LRT-set}. At first, it
may visually seem like the LRT-based set uses the threshold
$c_{\alpha,d}$, while the universal set uses $2\log(1/\alpha)$ which
is much smaller in high dimensions.  However, a key point to keep in
mind is that comparing the
numerators of the test statistics in both cases, the classical likelihood ratio set 
uses an \emph{in-sample
likelihood} and the split LRS confidence set 
uses an \emph{out-of-sample likelihood}. Hence, simply
comparing the thresholds does not suffice to draw a conclusion about the relative sizes of the confidence sets.
We next check that for
regular models, the size of the universal set indeed shrinks at the
right rate.



\section{Sanity Check: Regular Models}
\label{section::sanity}

Although universal methods are not needed
for well-behaved models, it is worth checking
their behavior in these cases.
We expect that $C_n$
would not have optimal size but we would hope that it still shrinks
at the optimal rate. 
We now confirm that this is true.

Throughout this example we treat the dimension as a fixed constant before subsequently turning our attention to 
an example where we more carefully track the dependence of the confidence set diameter on dimension. In this and subsequent sections
we use standard stochastic order notation for convergence in probability $o_p$, and boundedness in probability $O_p$ \cite{van2000asymptotic}. 
We make the following regularity assumptions (see for instance \cite{van2000asymptotic} for a detailed discussion of these conditions):
\begin{enumerate}
\item The statistical model is identifiable, i.e. for any $\theta \neq \theta^*$ it is the case that $P_{\theta} \neq P_{\theta^*}$.
The statistical model is differentiable in quadratic mean (DQM) at $\theta^*$, i.e. there exists a function $s_{\theta^*}$ such that:
\begin{align*}
\int & \left[\sqrt{p_\theta} - \sqrt{p_{\theta^*}} - \frac{1}{2} (\theta-\theta^*)^T s_{\theta^*}\sqrt{p_{\theta^*}}\right]^2 d\mu = \\
&o(\|\theta-\theta^*\|^2), \text{ as $\theta \rightarrow \theta^*$.}
\end{align*}
\item The parameter space $\Theta \subset \mathbb{R}^d$ is compact, and the log-likelihood is a smooth function
of $\theta$, i.e. there is a measurable function 
$\ell$ with $\sup_{\theta} P_{\theta} \ell^2 < \infty$ such that
for any $\theta_1, \theta_2 \in \Theta$:
\begin{align*}
|\log p_{\theta_1}(x) - \log p_{\theta_2}(x)| \leq \ell(x) \|\theta_1 - \theta_2\|. 
\end{align*}
\item A consequence of the DQM condition is that the Fisher information matrix 
\begin{align*}
I(\theta^*) := \mathbb{E}_{\theta^*} [s_{\theta^*} s_{\theta^*}^T],
\end{align*}
is well-defined, and we assume it is non-degenerate.
\end{enumerate}
Under these conditions the optimal confidence set has (expected) diameter
$O(1/\sqrt{n})$. Our first result shows that the same is true of the universal set, provided that the 
initial estimate $\widehat{\theta}_1$ is $\sqrt{n}$-consistent, i.e. $\|\widehat{\theta}_1 - \theta^*\| = O_p(1/\sqrt{n})$. 
Under the conditions of our theorem, this consistency 
condition is satisfied when $\widehat{\theta}_1$ is the MLE but our result
is more generally applicable.
\begin{theorem}
\label{thm:main}
Suppose that $\widehat{\theta}_1$ is a $\sqrt{n}$-consistent estimator of $\theta^*$.
Under the assumptions above, the split LRT confidence set has diameter $O_p(\sqrt{\log(1/\alpha)/n})$.
\end{theorem}
A proof of this result is in the supplement. At a high level, in order to bound the diameter of the split
LRT set it suffices to show that for any $\theta$ sufficiently far from $\theta^*$, it is the case that
\begin{align*}
\frac{\mathcal{L}_0(\theta)}{\mathcal{L}_0(\widehat{\theta}_1)} \leq \alpha. 
\end{align*}
In order to establish this, note that we can write this condition as:
\begin{align*}
\log \frac{\mathcal{L}_0(\theta)}{\mathcal{L}_0(\theta^*)} + \log \frac{\mathcal{L}_0(\theta^*)}{\mathcal{L}_0(\widehat{\theta_1})} \leq \log(\alpha). 
\end{align*}
Bounding first term requires showing if we consider any $\theta$ sufficiently far from $\theta^*$ its likelihood is small relative to the likelihood of $\theta^*$. We build on the 
work of Wong and Shen \cite{wong1995} who provide uniform upper bounds on the likelihood ratio under technical conditions which ensure that the statistical model 
is not too big. Conversely, to bound the second term we need to argue that if $\widehat{\theta}_1$ is sufficiently close to $\theta^*$, then it must be the case that their likelihoods 
cannot be too different. This in turn follows by exploiting the DQM condition.

\vspace{-0.05in}
\paragraph{Analyzing the non-parametric split LRT. } 
While our previous result focused on the diameter of the split LRT set in parametric problems, similar techniques also yield results in the non-parametric case. 
In this case, since we have no 
underlying parameter space 
it will be natural to measure the diameter of our confidence set in terms of some metric on probability distributions. We consider bounding the diameter of our
confidence set in the Hellinger metric. Formally, for two distributions $P$ and $Q$ the (squared) Hellinger distance is defined as:
\begin{align*}
H^2(P,Q) = \frac{1}{2} \int (\sqrt{dP} - \sqrt{dQ})^2. 
\end{align*}
We will also require the use of the $\chi^2$-divergence given by:
\begin{align*}
\chi^2(P,Q) = \int \Big( \frac{dP}{dQ} - 1\Big)^2 dQ,
\end{align*}
assuming that $P$ is absolutely continuous with respect to $Q$. Roughly, and analogous to our development in the parametric case,
in order to bound the diameter of the split LRT
confidence set, we need to ensure that our statistical model $\mathcal{P}$ is not too large, and further that our initial estimate $\widehat{p}_1$ is 
sufficiently close to $p^*$.
 
To measure the size of $\mathcal{P}$ we use its Hellinger bracketing entropy. 
Denote by $\log N(u, \mathcal{F})$ the Hellinger bracketing
entropy of the class of distributions $\mathcal{F}$ where the bracketing functions are separated by at most $u$ in the Hellinger distance (we refer to \cite{wong1995} for a precise definition).
We suppose that the bracketing entropy of $\mathcal{P}$ is not too large, i.e. for some $\epsilon_n > 0$ we have that for some constant $c > 0$,
\begin{align}
\label{eqn:con1}
\int_{\epsilon_n^2}^{\epsilon_n} \sqrt{ \log(N(u, \mathcal{P}))} du \leq c \sqrt{n} \epsilon_n^2.
\end{align}
Although we do not explore this in detail, we note in passing that the smallest value $\epsilon_n$ for which the above condition is satisfied 
provides an upper bound on the rate of convergence of the non-parametric MLE in the Hellinger distance \cite{wong1995}.
To characterize the quality of $\widehat{p}_1$ we use the $\chi^2$ divergence. Concretely, we suppose that:
\begin{align}
\label{eqn:con2}
\chi^2(p^*,\widehat{p}_1) \leq O_p(\eta_n^2).
\end{align}

\begin{theorem}
\label{thm:nonpar}
Under conditions~\eqref{eqn:con1} and~\eqref{eqn:con2}, 
the split LRT confidence set has Hellinger diameter upper bounded by $O_p(\eta_n + \epsilon_n + \sqrt{\log(1/\alpha)/n})$.
\end{theorem}

 \vspace{-0.1cm}

\paragraph{Comparing LRT to split LRT for the multivariate Normal case.}
In the previous calculation we treated the dimension of the parameter space as fixed. To understand 
the behavior of the method
as a function of dimension in the regular case,
suppose that
$Y_1,\ldots, Y_n \sim N_d(\theta,I)$
where $\theta\in\mathbb{R}^d$.
Recalling that we use $c_{\alpha,d}$ and $z_\alpha$ to denote the upper $\alpha$ quantiles of the $\chi^2_d$ and standard Gaussian respectively,
the usual confidence set for $\theta$ based on the LRT is
\begin{align*}
A_n &= 
\left\{\theta:\ \|\theta-\overline{Y}\|^2 \leq \frac{c_{\alpha,d}}{n}\right\}\\
&=
\left\{\theta:\ \|\theta-\overline{Y}\|^2 \leq \frac{d + \sqrt{2d} z_\alpha + o(\sqrt{d})}{n}\right\},
\end{align*}
where the second form follows from the Normal approximation of the $\chi^2_d$ distribution. 
For the Universal set, we use the sample average from $D_1$ as our initial estimate 
$\widehat\theta_1$.
Denoting the sample means $\overline{Y}_1$ and $\overline{Y}_0$ 
we see that 
$$
C_n = \left\{ \theta:\ \log {\cal L}_0(\overline{Y}_1) - \log {\cal L}_0(\theta) \leq \log (1/\alpha)\right\},
$$
which is the set of $\theta$ such that
$$
- \left(\frac{n}{2}\right) \frac{\|\overline{Y}_0 - \overline{Y}_1\|^2}{2} + 
\left(\frac{n}{2}\right) \frac{\|\theta-\overline{Y}_0\|^2}{2} \leq \log\left(\frac{1}{\alpha}\right).
$$
In other words, we may rewrite
$$
C_n = \left\{\theta:\ 
\|\theta-\overline{Y}_0\|^2 \leq \frac{4}{n}\log \left(\frac{1}{\alpha}\right) + \|\overline{Y}_0 - \overline{Y}_1\|^2 \right\}.
$$
Next, note that
$\|\overline{Y}_0 - \overline{Y}_1\|^2 = O_p(d/n)$,
so both sets have radii $O_p(d/n)$.
Precisely, the squared radius $R_n^2$ of $C_n$ is
\begin{align*}
R_n^2 &\stackrel{d}{=} ~
\frac{ 4 \log (1/\alpha) + 4 \chi_d^2}{n}\\
& \stackrel{d}{=} ~
\frac{ 4 \log (1/\alpha) + 4 d + \sqrt{32 d}\, Z + O_p(\sqrt{d})}{n},
\end{align*}
where $Z$ is an independent standard Gaussian.  So both their squared
radii share the same scaling with $d$ and $n$, and for large $d$
and constant $\alpha$, the squared radius of $C_n$ is about 4 times larger than that of $A_n.$

\section{Examples}
\label{sec:examples}

\paragraph{Mixture models.}
As a proof-of-concept, we do a small
simulation to check the type I error and power
for mixture models.
Specifically,
let
$Y_1,\ldots, Y_{2n} \sim P$ 
where $Y_i\in\mathbb{R}$.
We want to distinguish the hypotheses in~\eqref{eqn:mixture}.
For this brief example, we take
$k_0=1$ and $k_1=2$.

Finding a test that provably controls the type I error at a given level
has been elusive.
A natural candidate is the likelihood ratio statistic
but, as mentioned earlier, this has an intractable
limiting distribution. To the best of our knowledge,
the only practical test for the above hypothesis with
a tractable limiting distribution is the EM test due to
\cite{chen2009hypothesis}.
This very clever test is similar to the likelihood ratio test except that
it includes some penalty terms
and requires the maximization of some of the parameters to be restricted.
However,
the test requires choosing some tuning parameters
and, more importantly, it is restricted to one-dimensional problems.
There is no known confidence set
for mixture problems with guaranteed coverage properties.
Another approach is based on the bootstrap
\citep{mclachlan1987bootstrapping}
but there is no proof
of the validity of the bootstrap for mixtures.

\begin{figure}
\begin{center}
\includegraphics[scale=.4]{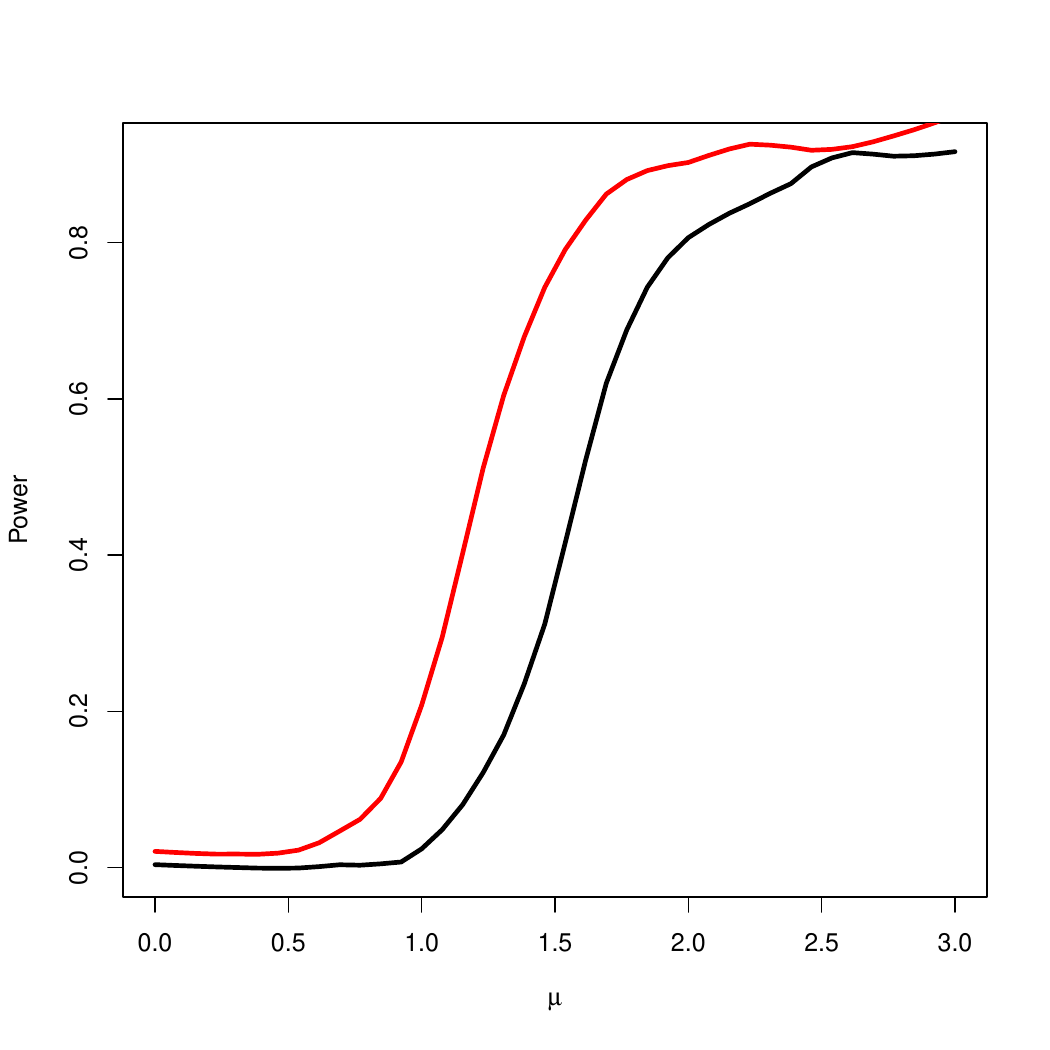}
\end{center}
\caption{\em The plot shows the power of the universal/bootstrap (black/red) tests for a simple Gaussian mixture, as the mean-separation $\mu$ varies ($\mu=0$ is the null).
The sample size is $n=200$ and the target level is $\alpha=0.1$.}
\label{fig::power}\vspace{-0.15in}
\end{figure}

Figure \ref{fig::power}
shows the power
of the test when $n=200$
and $\widehat\theta_1$ is the MLE under the full model
${\cal M}_2$.
The true model is taken to be
$(1/2) \phi(y; -\mu,1) + (1/2) \phi(y; \mu,1)$
where $\phi$ is a Normal density with mean $\mu$ and variance 1.
The null corresponds to $\mu=0$.
We take $\alpha = 0.1$ and
the MLE is obtained by the EM algorithm, which we assume converges on this simple problem.
Understanding the local and global convergence (and non-convergence) of the EM algorithm to the MLE 
is an active research area but is beyond the scope of this paper \cite[and references therein]{balakrishnan2017statistical,xu2016global,jin2016local}.
As expected, the test is conservative
with type I error near 0 but has reasonable power when $\mu > 1$.

Figure \ref{fig::power} also shows the power of the
bootstrap test \cite{mclachlan1987bootstrapping}.
Here, the $p$-value is obtained by bootstrapping the
LRS under the estimated null distribution.
As expected, this has higher power than the universal test
since it does not split the data. In this simulation, both tests control the type I error, but unfortunately the bootstrap test does not have any guarantee on the type I error, even asymptotically.
The lower power of the universal test is the price paid for having a finite-sample guarantee.
It is also worth noting that the bootstrap test requires running the EM algorithm for each bootstrap sample while the universal test
only requires one EM run.

\vspace{-0.05in}

\paragraph{Model selection using sieves.}
Sieves are a general approach to nonparametric inference.
A sieve \cite{shen1994convergence}
is a sequence of nested models
${\cal P}_1 \subset {\cal P}_2 \subset \cdots$.
If we assume that the true density $p^*$ is in ${\cal P}_j$ for some (unknown) $j$
then universal testing can be used to choose the model.
One possibility is to test
$H_j: p^* \in {\cal P}_j$ one by one for $j=1,2,\ldots$.
We reject $H_j$ if
\[
\prod_{i\in D_0} \frac{\widehat p_{j+1}(Y_i)}{\widehat p_{j}(Y_i)} > 1/\alpha,
\]
where $\widehat p_j$ is the MLE in model ${\cal P}_j$.
Then we take $\widehat j$ to be the first $j$
such that $H_j$ is not rejected, and proclaim that $p^* \in {\cal P}_{j}$ for some $j \geq \widehat j$.
Even though we test multiple different hypotheses and stop at a random $\widehat j$, this procedure still controls the type I error, meaning that
\[
P_{p^*}(p^* \in {\cal P}_{\widehat j - 1}) \leq \alpha,
\]
meaning that our proclamation is correct with high probability.  The
reason we do not need to correct for multiple testing is because a
type I error can occur only once we have reached the first $j$ such
that $p^*\in {\cal P}_j$.

A simple application is to choose the number of mixture components in a mixture model, as discussed in the previous example. Here are some other interesting examples in which the aforementioned ideas yield valid tests and model selection using sieves: 
\begin{enumerate}
\item[(a)] testing the number of hidden states in a hidden markov model (the MLE is computable using the Baum-Welch algorithm), 
\item[(b)] testing the number of latent factors in a factor model, 
\item[(c)] testing the sparsity level in a high-dimensional linear model $Y=X\beta+\epsilon$ (under $H_0: \beta$ is $k$-sparse, the MLE corresponds to best-subset selection). 
\end{enumerate}
Whenever we can compute the MLE (specifically, the likelihood it achieves), then we can run our universal test, and we can do model selection using sieves. We will later see that an upper bound of the maximum likelihood suffices, and is sometimes achievable by minimizing convex relaxations of the negative log-likelihood.


\paragraph{Nonparametric example: Shape constrained inference.}
A density $p$ is log-convave if
$p= e^g$ for some concave function $g$.
Consider testing
$H_0: $ $p$ is log-concave versus
$H_1: $ $p$ is not log-concave.
Let ${\cal P}_0$ be the set of log-concave densities and
let $\widehat p_0$ denote the nonparametric maximum
likelihood estimator over ${\cal P}_0$ computed using $D_0$
\cite{cule2010maximum}
which can be computed in polynomial time \cite{axelrod2019polynomial}.
Let $\widehat p_1$ be any nonparametric density estimator such
as the kernel density estimator
\cite{silverman2018density} fit on $D_1$.
In this case, the universal test is to reject $H_0$ when
$$
\prod_{i\in D_0} \frac{\widehat p_1(Y_i)}{\widehat p_0(Y_i)} > \frac{1}{\alpha}.
$$
To the best of our knowledge this is the first test for this problem with finite-sample guarantee.
Under the assumption that $p\in {\cal P}_0$,
the universal confidence set is
$$
C_n = \Biggl\{ p\in {\cal P}_0:\ \prod_{i\in D_0} p(Y_i) \geq \alpha \prod_{i\in D_0}\widehat p_1(Y_i) \Biggr\}.
$$
While the aforementioned test can be efficiently performed, the set $C_n$ may be hard to explicitly represent, but 
we can check if a distribution $p \in C_n$ efficiently.



\paragraph{Positive dependence (MTP$_2$).} The split LRT solves a variety of open problems related to testing for a general notion of positive dependence called  multivariate total positivity of order 2 or MTP$_2$ \cite{karlin1980classes}. The convex optimization problem of maximum likelihood estimation in Gaussian models under total positivity was recently solved  \cite{LUZ}, but in Example 5.8 and the following discussion, they state that the testing problem is still open. Given data from a multivariate distribution $p$, consider testing $H_0: p$ is Gaussian MTP$_2$, against $H_1: p$ is Gaussian (or an even more general alternative). Since Proposition 2.2 in their paper shows that the MLE under the null can be efficiently calculated, our universal test is applicable. 

In fact, calculating the MLE in any MTP$_2$ exponential family is a convex optimization problem \cite[Theorem 3.1]{lauritzen2019total}, thus making a test immediately feasible. As a particularly interesting special case, their Section 5.1 provides an algorithm for computing the MLE for MTP$_2$ Ising models. Testing $H_0: p$ is Ising MTP$_2$ against $H_1: p$ is Ising, is stated as an open problem in their Section 6, and is solved by our universal test. (We remark that even though the MTP$_2$ MLE is efficiently computable, evaluating the maximum likelihood in the Ising case may still take $O(2^d)$ time for a $d$-dimensional problem.)

Finally, MTP$_2$ can be combined with log-concavity, uniting shape constraints and dependence. General existence and uniqueness properties of the MLE for totally positive log-concave densities have been recently derived \cite{robeva2018maximum}, along with efficient algorithms to compute the MLE. Our methods immediately yield a test for $H_0: p$ is MTP$_2$ log-concave against $H_1: p$ is log-concave.

All the above models were singular, and hence the LRS has been hard to study. In some cases, its asymptotic null distribution is known to be a weighted sum of chi-squared distributions, where the weights are rather complicated properties of the distributions (usually unknown to the practitioner). In contrast, the split LRT is applicable without assumptions, and its validity is nonasymptotic.


\paragraph{Independence versus conditional independence.}
Consider data that are trivariate vectors of the form
$(X_{1i},X_{2i},X_{3i})$
which are modelled as trivariate Normal.
The goal is to test
$H_0: $ $X_1$ and $X_2$ are independent
versus
$H_1: $ $X_1$ and $X_2$ are independent given $X_3$.
The motivation for this test
is that this problem arises in the construction of
causal graphs.
It is surprisingly difficult to test these non-nested
hypotheses.
Indeed,
\cite{2019arXiv190601850G}
study carefully the subtleties of the problem
and they show that the limiting distribution of the LRS
is complicated and cannot be used for testing.
They propose a new test based on a concept called envelope distributions.
Despite the fact that the hypotheses are non-nested,
the universal test is applicable and can be used quite
easily for this problem.
Further, one can also flip $H_0$ and $H_1$ and test for conditional independence
in the Gaussian setting
 as well.
We leave it to future work to compare the power of the
universal test and the envelope test.


\paragraph{Crossfitting can beat splitting: uniform distribution.}  
In all previous examples, the split LRT is a reasonable choice.
However, in this example, the crossfit approach easily dominates the split approach. 
Note that this is a case where we would not recommend our universal tests
since there are well-studied standard confidence intervals in this model. 
The example is just meant to bring out the difference between
the split and crossfit approaches.

Suppose that $p_\theta$ is the uniform
density on $[0,\theta]$.  Let us take $\widehat\theta_1$ to be the MLE
from $D_1$.  Thus, $\widehat\theta_1$ is the maximum of the data points in
$D_1$.  Now ${\cal L}_0(\theta) = \theta^{-n} I(\theta\geq
\widehat\theta_0)$ where $\widehat\theta_0$ is the maximum of the data points
in $D_0$.  It follows that $C_n = [0,\infty)$ whenever $\widehat\theta_1 <
\widehat\theta_0$ which happens with probability 1/2.  The set $C_n$ has
the required coverage but is too large to be useful.  This happens
because the densities have different support.  
A similar phenomenon occurs when testing
$H_0: \theta \leq A$ versus $H_1: \theta \in \mathbb{R}^+$ for some fixed $A>0$, but not when
testing against $H_1: \theta > A$.
One can partially
avoid this behavior by choosing $\widehat \theta_1$ to not be the MLE.
However, the simplest way to avoid the degeneracy is to use the crossfit approach,
where we swap the roles of $D_0$ and $D_1$, and average the resulting test statistics.
Exactly one of two test statistics will be 0, and hence the average will be nonzero.
Further, it is easy to show that this test and resulting interval are rate-optimal,
losing a constant factor due to data splitting over the standard tests and interval constructions.
In more detail, the classical (exact) pivotal $1-\alpha$ confidence interval for $\theta$ is 
$C'_{2n} = [\widehat{\theta}, \widehat{\theta} (1/\alpha)^{1/(2n)}],$ where $\widehat{\theta}$ is
the maximum of all the data points. On the other hand, for $\widehat{\theta}_1, \widehat{\theta}_0$ defined above, 
assuming without loss of generality that $\widehat{\theta}_0 \leq \widehat{\theta}_1$ a direct calculation
shows that the crossfit interval
takes the form $C_n = [\widehat{\theta}_0, \widehat{\theta}_1 (2/\alpha)^{1/n}].$ Ignoring constants, both these 
intervals have expected length~$O(\theta \log(1/\alpha) /n)$. 

\vspace{-0.1in}
\section{De-randomization}

The universal method involves
randomly splitting the data
and the final inferences will depend
on the randomness of the split.
This may lead to instability, where different
random splits produce different results; in a related context,
this has been called the ``p-value lottery'' \citep{meinshausen2009p}.

We can get rid of or reduce the variability of our inferences,
at the cost of more computation
by using many splits, while maintaining validity of the method.
The key property that we used in both the universal confidence set and the split LRT is that
$\mathbb{E}_{\theta^*}[T_n]\leq 1$ where
$T_n= {\cal L}_0 (\widehat{\theta}_1)/{\cal L}_0 (\widehat{\theta})$.
Imagine that we obtained $B$ such statistics $T_{n,1}\ldots, T_{n,B}$ with the same property.
Let 
\begin{align*}
\overline{T}_n = B^{-1}\sum_{j=1}^B T_{n,j}.
\end{align*}
Then we still have that
$\mathbb{E}_{\theta^*} [\overline{T}_n] \leq 1$
and so inference using our universal methods can proceed using the combined statistic $\overline{T}_n$. Note that this is true regardless of the dependence between the statistics.

Using the aforementioned idea we can immediately design natural variants of the universal method:
\begin{itemize}
\item {\bf K-fold.} We can split the data once into $2 \leq K \leq n$ folds. Then repeat the following $K$ times: use $K-1$ folds to calculate $\widehat \theta_1$, and evaluate the likelihood ratio on the last fold. Finally, average the $K$ statistics. Alternatively, we could use one fold to calculate $\widehat \theta_1$ and evaluate the likelihood on the other $K-1$ folds.
\item {\bf Subsampling.} We do not need to split the data just once into $K$ folds. We can repeat the previous procedure for repeated random splits of the data into $K$ folds. We expect this to reduce variance that arises from the algorithmic randomness. 
\item {\bf All-splits.} We can remove all algorithmic randomness by considering all possible splits. While this is computationally infeasible, the potential statistical gains are worth studying. 
\end{itemize}
We remark that all these variants allow a large amount of flexibility. For example, in cross-fitting, $\widehat \theta_1$ need not be used the same way in both splits: it could be the MLE on one split, but a Bayesian estimator on another split. This flexibility could be useful if the user does not know which variant would lead to higher power in advance and would like to hedge across multiple natural choices. Similarly, in the $K$-fold version, if a user is confused whether to evaluate the likelihood ratio on one fold or on $K-1$ folds, then they can do both and average the statistics. 

Of course, with such flexibility comes the risk of an analyst cherry-picking the variant used after looking at the which form of averaging results in the highest LR (this would correspond to taking the maximum instead of the average of multiple variants), but this is a broader issue. For this reason (and this reason alone), the cross-fitting LRT proposed initially may be a useful default in practice, since it is both conceptually and computationally simple. We have already seen that (two-fold) cross-fit  inference improves over split inference drastically in the case of the uniform distribution discussed in the previous section. We leave a more detailed theoretical and empirical analysis of the power of these variants to future work. 


\section{Extensions}
\label{sec::ext}

\paragraph{Profile likelihood and nuisance parameters.} Suppose that we are interested in some function
$\psi = g(\theta)$.
Let 
\[
B_n = \{ \psi:\ C_n \bigcap g^{-1}(\psi) \neq \emptyset\},
\]
where we define $g^{-1}(\psi) = \{\theta: g(\theta) = \psi\}.$
By construction,
$B_n$ is a $1-\alpha$ confidence set for $\psi$.
Defining the profile likelihood function
\begin{equation}
{\cal L}_0^\dagger (\psi) = \sup_{\theta:~ g(\theta) = \psi}{\cal L}_0(\theta),
\end{equation}
we can rewrite $B_n$ as
\begin{equation}
B_n= \left\{ \psi:\ \frac{{\cal L}_0(\widehat\theta_1)}{{\cal L}_0^\dagger (\psi)} \leq \frac1\alpha \right\}.
\end{equation}
In other words, the same data splitting idea works for the profile likelihood too. As a particularly useful example, suppose $\theta=(\theta_u,\theta_n)$ where $\theta_n$ is a nuisance component, then we can define $g(\theta)=\theta_u$ to obtain a universal confidence set for only the component $\theta_u$ we care about.

\paragraph{Upper bounding the null maximum likelihood.} Computing the MLE and/or the maximum likelihood (under the null) is sometimes computationally hard.
Suppose one could come up with a relaxation $F_0$ of the null likelihood ${\cal L}_0$.  
This should be a proper relaxation in the sense that
\[
\max_\theta F_0(\theta) \geq \max_\theta {\cal L}_0(\theta).
\]
For example, ${\cal L}_0$ may be defined as $-\infty$ outside its domain, but $F_0$ could extend the domain.
As another example, instead of minimizing the negative log-likelihood which could be nonconvex and hence hard to minimize, 
we could minimize a convex relaxation.
In such settings, define 
\[
\widehat\theta^F_0 := \argmax_\theta F_0(\theta).
\]
If we define the test statistic
\[
T_n' := \frac{{\cal L}_0(\widehat \theta_1)}{F_0(\widehat \theta_0^F)},
\]
then the split LRT may proceed using $T_n'$ instead of $T_n$. This is because $F_0(\widehat \theta_0^F) \geq  {\cal L}_0(\widehat \theta_0)$, and hence $T_n' \leq T_n$.

One particular case when this would be useful is the following. While discussing sieves, we had mentioned that testing the sparsity level
in a high-dimensional linear model involves solving the best subset selection problem, which is NP-hard in the worst case. 
There exist well-known quadratic programming relaxations that are more computationally tractable. Another example is testing if a random graph is a stochastic block model, for which semidefinite relaxations of the MLE are well studied \cite{amini2018semidefinite}; similar situations arise in communication theory \cite{dahl2003approximate} and angular synchronization \cite{bandeira2017tightness}. 

The takeaway message is that it suffices to upper bound the maximum likelihood in order to perform inference.

\paragraph{Robustness via powered likelihoods.} It has been suggested by some authors 
\cite{royall2003interpreting,grunwald2012safe,holmes2017assigning,grunwald2017inconsistency,miller2019robust} 
that inferences can be made robust by replacing the likelihood ${\cal L}$
with the power likelihood ${\cal L}^\eta$ for some $0 < \eta < 1$.
Note that
$$
\mathbb{E}_\theta
\Biggl[
\Biggl(
\frac{ {\cal L}_0(\widehat\theta_1)}
{ {\cal L}_0(\theta)}
\Biggr)^\eta \, \Biggm| \, D_1\Biggr]=
\prod_{i\in D_0}\int
p_{\widehat\theta_1}^\eta (y_i)p_{\theta}^{1-\eta} (y_i) d y_i \leq 1,
$$
and hence all the aforementioned methods can be used with the robustified likelihood as well. (The last inequality follows because the $\eta$-Renyi divergence is nonnegative.)

\paragraph{Smoothed likelihoods.} Sometimes the MLE is not consistent or it may not exist since the likelihood function is unbounded, and a (doubly) smoothed likelihood has been proposed as an alternative \cite{seo2013universally}. For simplicity, consider a kernel $k(x,y)$ such that $\int k(x,y) dy = 1$ for any $x$; for example a Gaussian or Laplace kernel. For any density $p_\theta$,
 let its smoothed version be denoted 
\[
\widetilde p_\theta(y) := \int k(x,y) p_\theta(x) dx,
\]
Note that $\widetilde p_\theta$ is also a probability density.
Denote the smoothed empirical density based on $D_0$ as 
\[
\widetilde p_n := \frac1{|D_0|}\sum_{i \in D_0} k(X_i, \cdot ).
\]
Define the smoothed maximum likelihood estimator as the Kullback-Leibler (KL) projection of $\widetilde p_n$ onto $
\{\widetilde p_\theta\}_{\theta \in \Theta_0}$:
\[
\widetilde \theta_0 := \arg\min_{\theta \in \Theta_0} K(\widetilde p_n, \widetilde p_\theta),
\]
where $K(P,Q)$ denotes the KL divergence between $P$ and $Q$.
If we define the smoothed likelihood on $D_0$ as
\[
\widetilde{\cal L}_0(\theta) := \prod_{i\in D_0} \exp \int k(X_i, y) \log \widetilde p_\theta(y) dy,
\]
then it can be checked that  $\widetilde \theta_0$ maximizes the smoothed likelihood, that is
$
\widetilde \theta_0 = \arg\max_{\theta \in \Theta_0} \widetilde{\cal L}_0(\theta).
$
As before, let $\widehat \theta_1 \in \Theta$ be any estimator based on $D_1$. The smoothed split LRT is defined analogous to \eqref{eq:splitLRT} as:
\begin{equation}\label{eq:smoothed-splitLRT}
\text{reject $H_0$ if } \widetilde U_n > 1/\alpha, \text{ where } \widetilde U_n = \frac{\widetilde{\cal L}_0(\widehat\theta_1)}{\widetilde{\cal L}_0(\widetilde\theta_{0})}.
\end{equation}
We now verify that the smoothed split LRT controls type-1 error. First, for any fixed $\psi \in \Theta$, we have
\begin{align*}
\mathbb{E}_{\theta^*}\left[\frac{\widetilde{\cal L}_0(\psi)}{\widetilde{\cal L}_0(\widetilde\theta_{0})}\right] 
&\stackrel{\text{(i)}}{\leq} \mathbb{E}_{\theta^*}\left[\frac{\widetilde{\cal L}_0(\psi)}{\widetilde{\cal L}_0(\theta^*)}\right] \\
& \hspace{-1cm}= \prod_{i \in D_0} 
\int \exp \left( \int k(x,y) \log \frac{\widetilde p_{\psi}(y)}{\widetilde p_{\theta^*}(y) } dy \right) p_{\theta^*}(x) dx\\
&\hspace{-1cm}\stackrel{\text{(ii)}}{\leq} 
\int \left( \int k(x,y) \frac{\widetilde p_{\psi}(y)}{\widetilde p_{\theta^*}(y) } dy \right) p_{\theta^*}(x) dx\\
&\hspace{-1cm}= \int \left(  \frac{\int k(x,y) p_{\theta^*}(x) dx }{\widetilde p_{\theta^*}(y) } \right) \widetilde p_{\psi}(y) dy \\
&\hspace{-1cm}= \int \widetilde p_{\psi}(y) dy = 1.
\end{align*}
Above, step $\text{(i)}$ is because $\widetilde \theta_0$ maximizes the smoothed likelihood, and step $\text{(ii)}$ follows by Jensen's inequality.
An argument mimicking equations \eqref{eq:unit-expectation} and \eqref{eq:markov} completes the proof. As a last remark, similar to the unsmoothed case, note that upper bounding the smoothed maximum likelihood under the null also suffices.

\vspace{-0.05in}
\paragraph{Conditional likelihood for non-i.i.d. data.} Our presentation so far has assumed that the data are drawn i.i.d. from some distribution under the null. However, this is not really required (even under the null), and was assumed for expositional simplicity. All that is needed is that we can calculate the likelihood on $D_0$ conditional on $D_1$ (or vice versa). For example, this could be tractable in models involving sampling without replacement from an urn with $M \gg n$ balls. Here $\theta$ could represent the unknown number of balls of different colors. Such hypergeometric sampling schemes result in non-i.i.d. data, but conditional on one subset of data (for example how many red, green and blue balls were sampled from the urn in that subset), one can evaluate the conditional likelihood of the second half of the data and maximize it, rendering it possible to apply our universal tests and confidence sets.

\vspace{-0.1in}
\section{Misspecification, and convex model classes}
\label{sec:mis}

There are some natural
examples of convex model classes \cite{li1999estimation,hoff2003nonparametric}, including (A) all mixtures (potentially infinite) of a set of base distributions, (B) distributions with the first moment specified/bounded and possibly other moments bounded (eg: first moment equals zero, second moment bounded by one), (C) the set of (coordinatewise) monotonic densities with the same support, (D) unimodal densities with the same mode, (E) densities that are symmetric about the same point, (F) distributions with the same median or multiple quantiles (eg: median equals zero, $0.9$-quantile equals two), (G) the set of all $K$-tuples $(P_1,\dots,P_K)$ of distributions satisfying a fixed partial stochastic ordering (eg: all triplets $(P_1,P_2,P_3)$ such that $P_1 \preceq P_2$ and $P_1 \preceq P_3$, where $\preceq$ is the usual stochastic ordering), (H) the set of convex densities with the same support. Some cases like (F) and (G) also result in weakly closed convex sets, as does case (B) for a specified mean. (Several of these examples also apply in discrete settings such as constrained multinomials.)

It is often possible to calculate the MLE over these convex model classes using convex optimization, for example see \cite{brunk1966maximum,dykstra1989nonparametric} for case (G). This renders our universal tests and confidence sets immediately applicable. However, in this special case, it is also possible to construct new tests, and the universal confidence set has some nontrivial guarantees if the model is misspecified.

\paragraph{Model misspecification.} 
Suppose the data come from a distribution $Q$ with density $q \notin \mathcal{P}_\Theta \equiv \{p_\theta\}_{\theta \in \Theta}$, meaning that the model is misspecified and the true distribution does not belong to the considered model. In this case, what does the universal set $C_n$ defined in \eqref{eqn:universal} contain?
We will answer this question when the set of measures/densities $\mathcal{P}_\Theta$ is convex.  Define the Kullback-Leibler divergence of $q$ from $\mathcal{P}_\Theta$ as 
\[
K(q,\mathcal{P}_\Theta) := \inf_{\theta \in \Theta} K(q, p_\theta).
\]
Following Definition~4.2 in Li's PhD thesis \cite{li1999estimation}, a function $p^* \equiv p^*_{q \to \Theta}$ is called the reversed information projection (RIPR) of $q$ onto $\mathcal{P}_\Theta$ if for every sequence $p_n$ with $K(q,p_n) \to K(q,\mathcal{P}_\Theta)$, we have $\log p_n \to \log p^*$ in $L^1(Q)$. Theorem~4.3 in \cite{li1999estimation} proves that $p^*$ exists and is unique, satisfies $K(q,p^*) =  K(q,\mathcal{P}_\Theta)$, and
\begin{equation}
\label{eq:ripr}
\forall \theta \in \Theta,~ \mathbb{E}_{Y \sim q}\left[\frac{p_\theta(Y)}{p^*(Y)} \right] \leq 1.
\end{equation}
The above statement can be loosely interpreted as ``if the data come from $q \notin  \mathcal{P}_\Theta$, its RIPR $p^*$ will have higher likelihood than any other model in expectation''. We discuss this condition further at the end of this subsection.

It might be reasonable to ask whether the universal set contains $p^*$.
For various technical reasons (detailed in \cite{li1999estimation}) 
it is not the case, in general, that $p^*$ belongs to the collection $\mathcal{P}_\Theta$.
Since the universal set only considers densities in $\mathcal{P}_\Theta$ by construction, it cannot possibly contain $p^*$ in general.
However, when $p^*$ is a density in $\mathcal{P}_\Theta$, then it is indeed covered by our universal set.
\begin{proposition}
Suppose that the data come from $q \notin \mathcal{P}_\Theta$. If $\mathcal{P}_\Theta$ is convex and there exists a density $p^* \in \mathcal{P}_\Theta$ such that $K(q,p^*) = \inf_{\theta \in \Theta} K(q, p_\theta)$, then we have $P_q(p^* \in C_n) \geq 1-\alpha$.
\end{proposition}

The proof is short. Examining the proof of Theorem~\ref{thm:universal-set}, we must simply verify that for each $i \in D_0$, we have
\[
\mathbb{E}_{q}\left[\frac{p_{\widehat \theta_1}(Y_i)}{p^*(Y_i)} \right] \leq 1,
\]
which follows from \eqref{eq:ripr}. 
Here is a heuristic argument for why \eqref{eq:ripr} holds when $p^* \in \mathcal{P}_\Theta$. For any $\theta \in \Theta$, note that $K(q,\mathcal{P}_\Theta) = K(q,p^*) = \min_{\alpha \in [0,1]} K(q, \alpha p^* + (1-\alpha)p_\theta)$ since $\mathcal{P}_\Theta$ is convex. The KKT condition for this optimization problem is that gradient with respect to $\alpha$ is negative at $\alpha=1$ (the minimizer). Exchanging derivative and integral immediately yields \eqref{eq:ripr}. This argument is fomalized in Chapter 4 of \cite{li1999estimation}.


\paragraph{An alternate split LRT (RIPR Split LRT).}
We return back to the well-specified case for the rest of this paper.
First note that the fact in \eqref{eq:ripr} can be rewritten as
\begin{equation}
\label{eq:ripr2}
\forall \theta \in \Theta,~  \mathbb{E}_{Y \sim p_\theta}\left[\frac{q(Y)}{p^*(Y)} \right] \leq 1,
\end{equation}
which is informally interpreted as ``if the data come from $p_\theta$, then any alternative $q \notin \mathcal{P}_\Theta$ will have lower likelihood than its RIPR $p^*$ in expectation''.
This motivates the development of an alternate \emph{RIPR split LRT} to test composite null hypotheses that is defined as follows. As before, we divide the data into two parts, $D_0$ and $D_1$, and let $\widehat \theta_1 \in \Theta_1$ be any estimator found using only $D_1$. Now, define $p_0^*$ to be the RIPR of $p_{\widehat \theta_1}$ onto the null set $\{p_\theta\}_{\theta \in \Theta_0}$.  The RIPR split LRT rejects the null if
\[
R_n \equiv \prod_{i \in D_0} \frac{p_{\widehat \theta_1}(Y_i)}{p_0^*(Y_i)} > 1/\alpha.
\]
The main difference from the original MLE split LRT, is that earlier we had ignored $\widehat \theta_1$, and simply calculated the MLE $\widehat \theta_0$ under the null based on $D_0$.
\begin{proposition}
If $\{p_\theta\}_{\theta \in \Theta}$ is a convex set of densities, then $\sup_{\theta_0 \in \Theta_0} P_{\theta_0}(R_n > 1/\alpha) \leq \alpha$.
\end{proposition}
The fact that $p_0^*$ is potentially not an element of $\{p_{\theta}\}_{\theta \in \Theta_0}$ does not matter here. The validity of the test follows exactly the same logic as the MLE split LRT, observing that~\eqref{eq:ripr2} implies that for any true $\theta^* \in \Theta_0$, we have
\[
\mathbb{E}_{p_{\theta^*}}\left[\frac{p_{\widehat \theta_1}(Y_i)}{p_0^*(Y_i)} \right] \leq 1.
\]
Without sample splitting and with a fixed alternative distribution, the RIPR LRT has been recently studied \cite{grunwald_safe_2019}. 
When $\mathcal{P}_\Theta$ is convex and the RIPR split LRT is implementable, meaning that it is computationally feasible to find the RIPR or evaluate its likelihood, then this test can be more powerful than the MLE split LRT. Specifically, if the RIPR is actually a density in the null set, then
\[
R_n = \prod_{i\in D_0} \frac{p_{\widehat \theta_1}(Y_i)}{p_0^*(Y_i)} ~\geq~ \prod_{i \in D_0} \frac{p_{\widehat \theta_1}(Y_i)}{p_{\widehat \theta_0}(Y_i)} = U_n,
\] since $\widehat \theta_0$ maximizes the denominator among null densities.
Because of the restriction to convex sets, and since there exist many more subroutines to calculate the MLE over a set than to find the RIPR, the MLE split LRT is more broadly applicable than the RIPR split LRT.

\vspace{-0.1in}
\section{Anytime $p$-values and confidence sequences}
\label{sec:seq}

Just like the sequential likelihood ratio test \cite{wald_sequential_1945} extends the LRT, the split LRT has a simple sequential extension. Similarly, the confidence set can be 
extended to a ``confidence sequence'' \cite{darling1967confidence}. 

Suppose the split LRT failed to reject the null. Then we are allowed to collect more data and update the test statistic (in a particular fashion), and check if the updated statistic crosses $1/\alpha$. If it does not, we can further collect more data and reupdate the statistic, and this process can be repeated indefinitely. Importantly we do not need any correction for repeated testing; this is primarily because the statistic is upper bounded by a nonnegative martingale. We describe the procedure next in the case when each additional dataset is of size one, but the same idea applies when we collect data in groups.



\paragraph{The running MLE sequential LRT.} Consider the following, more standard, sequential testing/estimation setup. We observe an i.i.d. sequence $Y_1,Y_2,\dots$ from $P_{\theta^*}$. We would like to test the hypothesis in \eqref{eq:test}. Let $\widehat \theta_{1,t-1}$ be any \emph{non-anticipating} estimator based on the first $t-1$ samples, for example the MLE, $\argmax_{\theta \in \Theta_1} \prod_{i=1}^{t-1} p_\theta(Y_i)$, or a regularized version of it to avoid misbehavior at small sample sizes. Denote the null MLE as
\[
\widehat \theta_{0,t} = \argmax_{\theta \in \Theta_0} \prod_{i=1}^t p_\theta(Y_i)
\] 
At any time $t$, reject the null and stop if
\[
M_t := \frac{\prod_{i=1}^t p_{\widehat \theta_{1,i-1}} (Y_i) }{ \prod_{i=1}^t p_{\widehat \theta_{0,t}}(Y_i) } > 1/\alpha.
\]
This test is computationally expensive: we must calculate $\widehat \theta_{1,t-1}$ and $\widehat \theta_{0,t}$ at each step. In some cases, these may be quick to calculate by warm-starting from $\widehat \theta_{1,t-2}$ and $\widehat \theta_{0,t-1}$.
For example, the updates can be done in constant time for exponential families, since the MLE is often a simple function of the sufficient statistics. However, even in these cases, the denominator takes time $O(t)$ to recompute at step $t$.

The following result shows that with probability at least $1 - \alpha$, this test will never stop under the null.  Let $\tau_\theta$ denote the stopping time when the data is drawn from $P_\theta$, which is finite only if we stop and reject the null.
\begin{theorem}
The running MLE LRT has type I error at most $\alpha$, meaning that $\sup_{\theta^* \in \Theta_0} P_{\theta^*}(\tau_{\theta^*} < \infty) \leq \alpha$.
\end{theorem} 
The proof involves the simple observation that under the null, $M_t$ is upper bounded by a nonnegative martingale $L_t$ with initial value one. Specifically, define the (oracle) process starting with $L_0 := 1$ and
\begin{equation}\label{eq:oracle-martingale}
L_t := \frac{\prod_{i=1}^t p_{\widehat \theta_{i-1}} (Y_i) }{ \prod_{i=1}^t p_{\theta^*}(Y_i) } ~\equiv~ L_{t-1} \frac{p_{\widehat \theta_{t-1}} (Y_t) }{ p_{\theta^*}(Y_t) }.
\end{equation}
Note that under the null, we have $M_t \leq L_t$ because $\widehat \theta_{0,t}$ and $\theta^*$ both belong to $\Theta_0$, but the former maximizes the null likelihood (denominator). Further, it is easy to verify that $L_t$ is a nonnegative martingale with respect to the natural filtration $\mathcal{F}_t = \sigma(Y_1,\dots,Y_t)$. Indeed,
\begin{align*}
\mathbb{E}_{\theta^*}[L_t | \mathcal{F}_{t-1}]  &=
\mathbb{E}_{\theta^*}\left[ \frac{\prod_{i=1}^t p_{\widehat \theta_{i-1}} (Y_i) }{ \prod_{i=1}^t p_{\theta^*}(Y_i) } 
\Biggm| \mathcal{F}_{t-1} \right]\\
&= L_{t-1} \mathbb{E}_{\theta^*}\left[ \frac{ p_{\widehat \theta_{t-1}} (Y_t) }{ p_{\theta^*}(Y_t) } 
\Biggm| \mathcal{F}_{t-1} \right] = L_{t-1},
\end{align*}
where the last equality mimics \eqref{eq:unit-expectation}. To complete the proof, we note that the type I error of the running MLE LRT is simply bounded as
\begin{align*}
P_{\theta^*}(\exists t \in \mathbb{N}: M_t > 1/\alpha) &\leq P_{\theta^*}(\exists t \in \mathbb{N}: L_t > 1/\alpha)\\
&\stackrel{\text{(i)}}{\leq} \mathbb{E}_{\theta^*}[L_0] \cdot \alpha ~ = ~\alpha,
\end{align*}
where step $\text{(i)}$ follows by Ville's inequality \cite{ville_etude_1939,howard_exponential_2018}, a time-uniform version of Markov's inequality for nonnegative supermartingales. 

Naturally, this test does not have to start at $t=1$ when only one
sample is available, meaning that we can set $M_0=M_1=\dots=M_{t_0}=1$
for the first $t_0$ steps and then begin the
updates. Similarly, $t$ need not represent the time at which the
$t$-th sample was observed, it can just represent the $t$-th
recalculation of the estimators (there may be multiple samples
observed between $t-1$ and $t$).

\paragraph{Anytime-valid $p$-values.} We can also get a $p$-value that is uniformly valid over time. Specifically, both $p_t = 1/M_t$ and $\bar p_t = \min_{s \leq t} 1/M_s$ may serve as $p$-values. 
\begin{theorem}
For any random time $T$, not necessarily a stopping time, $\sup_{\theta^* \in \Theta_0} P_{\theta^*}(\bar p_T \leq x) \leq x$ for $x \in [0,1]$.
\end{theorem}
The aforementioned property is equivalent to the statement that under the null $P(\exists t \in \mathbb{N}: \bar p_t \leq \alpha) \leq \alpha$, and its proof follows by substitution immediately from the previous argument.
Naturally $\bar p_t \leq p_t$, but from the perspective of designing a level $\alpha$ test they are equivalent, because the first time that $p_t$ falls below $\alpha$ is also the first time that $\bar p_t$ falls below $\alpha$. 
The term ``anytime-valid'' is used because, unlike typical $p$-values, these are valid at (data-dependent) stopping times, or even random times chosen post-hoc. Hence, inference is robust to ``peeking'', optional stopping, and optional continuation of experiments. 
Such anytime $p$-values can be inverted to yield confidence sequences, as described below.

\paragraph{Confidence sequences.} A confidence sequence for $\theta^*$ is an infinite sequence of confidence intervals that are all simultaneously valid. Such confidence intervals are valid at arbitrary stopping times, and also at other random data-dependent times that are chosen post-hoc. In the same setup as above, but without requiring a null set $\Theta_0$, define the running MLE likelihood ratio process
\[
R_t(\theta) := \frac{\prod_{i=1}^t p_{\widehat \theta_{1,i-1}} (Y_i) }{ \prod_{i=1}^t p_{\theta}(Y_i) }.
\]
Then, a confidence sequence for $\theta^*$ is given by
\[
C_t := \{\theta: R_t(\theta) \leq 1/\alpha \}.
\]
In fact, the running intersection $\bar C_t = \bigcap_{s \leq t} C_s$ is also a confidence sequence; note that $\bar C_t \subseteq C_t$.

\begin{theorem}
$C_t$ and $\bar C_t$ are confidence sequences for $\theta^*$, meaning that $P_{\theta^*}(\exists t \in \mathbb{N}: \theta^* \notin \bar C_t) \leq \alpha$. Equivalently, $P_{\theta^*}(\theta^* \in C_\tau) \geq 1-\alpha$ for any stopping time $\tau$, and also $P_{\theta^*}(\theta^* \in C_T) \geq 1-\alpha$ for any arbitrary random time $T$.
\end{theorem}

The proof is straightforward. First, note that $\theta^* \notin \bar C_t$ for some $t$ if and only if $\theta^* \notin C_t$ for some $t$. Hence,
\[
P_{\theta^*}(\exists t\in \mathbb{N}: \theta^* \notin C_t) = P_{\theta^*}(\exists t\in \mathbb{N}: R_t(\theta^*) > 1/\alpha) \leq \alpha,
\]
where the last step uses, as before, Ville's inequality for the martingale $R_t(\theta^*) \equiv L_t$ from \eqref{eq:oracle-martingale}. The fact that the other two statements in the theorem are equivalent to the first follows from recent work \cite{Aadi}.

\paragraph{Duality.} It is worth remarking that confidence sequences are dual to anytime $p$-values, just like confidence intervals are dual to standard $p$-values, in the sense that a $(1-\alpha)$ confidence sequence can be formed by inverting a family of level $\alpha$ sequential tests (each testing a different point in the space), and a level $\alpha$ sequential test for  a composite null set $\Theta_0$ can be obtained by checking if the $(1-\alpha)$ confidence sequence intersects the null set $\Theta_0$. 

In fact, our constructions of $p_t$ and $C_t$ (without running minimum/intersection) obey the same property: $p_t < \alpha$ only if $C_t \cap \Theta_0 = \emptyset$, and the reverse implication follows if $\Theta_0$ is closed. To see the forward implication, assume that there exists some element $\theta' \in C_t \cap \Theta_0$. Since $\theta' \in C_t$, we have $R_t(\theta') \leq 1/\alpha$. Since $\theta' \in \Theta_0$, we have $\inf_{\theta^* \in \Theta_0} R_t(\theta^*) \leq 1/\alpha$. This last condition can be restated as $M_t \leq 1/\alpha$, which means that $p_t \geq \alpha$. 

It is also possible to obtain an anytime $p$-value from a family of confidence sequences at different $\alpha$, by defining $p_t$ as the smallest $\alpha$ for which $C_t \equiv C_t(\alpha)$ intersects $\Theta_0$.

\paragraph{Extensions.} All the extensions from Section~\ref{sec::ext} extend immediately to the sequential setting. One can handle nuisance parameters using profile likelihoods; this for example leads to sequential $t$-tests (for the Gaussian family, with the variance as a nuisance parameter), which also yield confidence sequences for the Gaussian mean with unknown variance. Non-i.i.d. data, such as in sampling without replacement, can be handled using conditional likelihoods, and robustness can be increased with powered likelihoods. In these situations, the corresponding underlying process $L_t$ may not be a martingale, but a supermartingale. Also, as before, we may also use upper bounds on the maximum likelihood at each step (perhaps minimizing convex relaxations of the negative log-likelihood), or smooth the likelihood if needed.

Such confidence sequences have been developed under very general
nonparametric, multivariate, matrix and continuous-time settings using
generalizations of the aforementioned supermartingale technique; see
\cite{howard_exponential_2018,Aadi,howard2019sequential}. The
connection between anytime-valid $p$-values, $e$-values, safe tests, peeking, confidence sequences, and the properties of optional stopping and continuation have been explored recently \cite{Aadi,johari_peeking_2017,grunwald_safe_2019,shafer_test_2011}. The connection to the present work is that when run sequentially, our universal (MLE or RIPR) split LRT yields an anytime-valid $p$-value, an $e$-value, a safe test, can be inverted to form universal confidence sequences, and are valid under optional stopping and continuation, and these are simply because the underlying process of interest is bounded by a nonnegative (super)martingale. This line of research began over 50 years
ago by Robbins, Darling, Lai and Siegmund
\citep{darling1967confidence,robbins_statistical_1970,robbins_class_1972,lai_confidence_1976,lai_boundary_1976}. 
In fact, for testing \emph{point} nulls, the running MLE (or non-anticipating) martingale was suggested in passing by Wald \cite[Eq. 10:10]{wald1947sequential}, analyzed in depth by \cite{robbins_class_1972,robbins_expected_1974} 
where connections were shown to the mixture sequential probability ratio test. 
These ideas have been utilized in changepoint detection for both point nulls  \cite{lorden2005nonanticipating} and composite nulls \citep{vexler2008martingale}.


\section{Conclusion}

Inference based on the split
likelihood ratio statistic (and variants)
leads to simple tests and confidence sets
with finite-sample guarantees.
Our methods are most useful
in problems where standard asymptotic methods
are difficult/impossible to apply, such as complex composite null testing problems or nonparametric confidence sets.
Going forward,
we intend to run simulations
in a variety of models
to study the power of the test
and the size of the confidence sets, and study their optimality
in special cases. 
We do not expect the test to be rate optimal in all cases,
but it might have analogous properties to the generalized LRT. 
It would also be interesting to extend these methods (like the profile likelihood variant) to semiparametric problems 
where there is a finite dimensional parameter of interest and an infinite dimensional nuisance parameter.


\noindent
\paragraph{Acknowledgments.} We thank Caroline Uhler and Arun K. Kuchibhotla
 for references to open problems in
shape constrained inference, Ryan Tibshirani for suggesting
the relaxed likelihood idea. We are grateful to Bin Yu, Hue Wang and Marco Molinaro for
helpful feedback which motivated parts of Section~\ref{sec:mis}.
Thanks to the reviewers and Dennis Boos for helpful suggestions, and to Len Stefanski
for pointing us to work on smoothed likelihoods.




\bibliographystyle{plainnat}
\bibliography{Universal}


%

\appendix

\section{Proof of Theorem~\ref{thm:main}}
Throughout our proof we will use $c_1, c_2, \ldots$ to denote positive universal constants which may change from line to line.
Define $d_n$ to be the diameter of the split-LRT set. Fix $\kappa > 0$, we want to show that, for some finite $M > 0$, 
\begin{align*}
P(d_n \geq M\sqrt{\log(1/\alpha)/n}) \leq \kappa,
\end{align*}
for all $n$ large enough. Equivalently, we 
want to show that for any $\theta$ such that $\|\theta - \theta^*\| \geq \frac{M}{2}\sqrt{\log(1/\alpha)/n}$ we have that:
\begin{align}
\label{eqn:main}
\frac{\mathcal{L}_0(\theta)}{\mathcal{L}_0(\widehat{\theta})} \leq \alpha, 
\end{align}
with probability at least $1 - \kappa$.
We know that $\|\widehat{\theta} - \theta^*\| = O_p(1/\sqrt{n})$, so let us consider the event where
$\|\widehat{\theta} - \theta^*\| \leq c_1/\sqrt{n}$ which happens with probability at least $1 - \kappa/3$ for sufficiently large $c_1 > 0$. We condition on this event throughout the remainder 
of our proof.

Now let us focus on showing~\eqref{eqn:main}.
This is equivalent to showing that,
\begin{align*}
\log \frac{\mathcal{L}_0(\theta)}{\mathcal{L}_0(\theta^*)} + \log \frac{\mathcal{L}_0(\theta^*)}{\mathcal{L}_0(\widehat{\theta})} \leq \log (\alpha). 
\end{align*}
The bulk of the technical analysis is in analyzing each of these terms. We show the following bounds:
\begin{lemma}
\label{lem:XX}
We have the following bounds:
\begin{enumerate}
\item There is some fixed constant $c_2 > 0$ such that for any $\epsilon \geq c_2/\sqrt{n}$, 
\begin{align}
\label{eqn:claim1}
\sup_{\theta: \|\theta - \theta^*\| \geq \epsilon} \log \frac{\mathcal{L}_0(\theta)}{\mathcal{L}_0(\theta^*)} \leq O_p( - n \epsilon^2).
\end{align}

\item Furthermore, if $\|\widehat{\theta} - \theta^*\| \leq c_1/\sqrt{n}$ for some fixed constant $c_1 > 0$ then,
\begin{align}
\label{eqn:claim2}
\log \frac{\mathcal{L}_0(\theta^*)}{\mathcal{L}_0(\widehat{\theta})} \leq O_p(1).
\end{align}

\end{enumerate}
\end{lemma}
With these results in place the remainder of the proof is straightforward. In particular, combining 
each of these convergence in probability results, together with a union bound we obtain that 
for  any $\theta$ such that $\|\theta - \theta^*\| \geq \frac{M}{2}\sqrt{\log(1/\alpha)/n}$, for a sufficiently large
constant $M > 0$, we have that with probability at least $1 - \kappa$,
\begin{align*}
\log \frac{\mathcal{L}_0(\theta)}{\mathcal{L}_0(\widehat{\theta})} \leq \log \alpha,
\end{align*}
as desired. To complete the proof it remains to prove Lemma~\ref{lem:XX}, and we prove each of its claims in turn.

\subsection{Proof of Claim~\eqref{eqn:claim1}}
In the proof of this result, it will be convenient to relate a natural metric on the underlying distributions (the Hellinger metric), and a natural 
metric on the underlying parameter space (the $\ell_2$ metric). We have the following result:
\begin{lemma}
\label{lem:hell_main}
Under the assumptions of our theorem:
\begin{enumerate}
\item There is a universal constant $c_1 > 0$ such that,
$H(p_{\theta_1},p_{\theta_2}) \leq c_1 \|\theta_1 - \theta_2\|$.

\item There are universal constants $c_1, c_2 > 0$ such that,
for any $\theta \in \Theta$ if $H(p_{\theta}, p_{\theta^*}) \leq c_1$, then
$H(p_{\theta}, p_{\theta^*}) \geq c_2 \|\theta - \theta^*\|.$
\end{enumerate}
\end{lemma}
Roughly, this result guarantees us that the Hellinger distance is always upper bounded by the $\ell_2$ distance, and further 
that in a small neighborhood of $\theta^*$ the Hellinger distance is also lower bounded by the $\ell_2$ distance.
We defer the proof of this result to Section~\ref{sec:proof_of_lemma_hell_main}, and now turn our attention to bounding the diameter of the split-LRT set.


We build on the results of Wong and Shen \cite{wong1995} who characterize the behaviour of likelihood ratios under assumptions on the Hellinger bracketing entropy 
of the underlying statistical model.
Towards this we first bound the local metric entropy of our statistical model in the following lemma. We denote by $\log N(u, \mathcal{F})$ the Hellinger bracketing
entropy of the class of distributions $\mathcal{F}$ where the bracketing functions are separated by at most $u$ in the Hellinger distance. We denote by $\mathcal{F}$ the collection of distributions $P_{\theta}$ for $\theta \in \Theta$.
\begin{lemma}
\label{lem:local_entropy}
There exist constants $c_1, c_2, c_3  > 0$ such that for any $s \geq c_1/\sqrt{n} > 0$, 
\begin{align*}
\int_{s^2/8}^{\sqrt{2}s} \sqrt{ \log N(u/c_2,\mathcal{F} \cap \left\{ H^2(p_{\theta},p_{\theta^*}) \leq s^2 \right\} )} du \leq c_3 \sqrt{n}s^2.
\end{align*}
\end{lemma}
With this local bracketing entropy bound in place Theorem 2 of Wong and Shen \cite{wong1995} yields the following conclusion: there exist constants $c_4, c_5, c_6 > 0$ such that for any $\epsilon \geq c_4/\sqrt{n}$, 
\begin{align*}
P_{\theta^*}\Big(\sup_{H(p_{\theta},p_{\theta^*})  \geq \epsilon} \log \left(\frac{\mathcal{L}_0(\theta)}{\mathcal{L}_0(\theta^*)}\right) \geq - c_5 n\epsilon^2 \Big) \leq 4 \exp (-c_6 n \epsilon^2).
\end{align*}
The desired claim follows immediately.
\subsubsection{Proof of Lemma~\ref{lem:hell_main}}
\label{sec:proof_of_lemma_hell_main}

\leavevmode \\
\noindent {\bf Proof of Claim 1: } We begin with our regularity condition: 
\begin{align*}
\left|\log \frac{p_{\theta_1}}{p_{\theta_2}}\right| \leq \ell(x) \|\theta_1 - \theta_2\|.
\end{align*}
Using the inequality that for $x \geq 0$, 
\begin{align*}
1 - 1/x \leq \log x,
\end{align*}
we obtain that,
\begin{align*}
1 - \sqrt{\frac{p_{\theta_2}}{p_{\theta_1} }} \leq  \frac{1}{2} \log \frac{p_{\theta_1}}{p_{\theta_2}}  \leq \frac{\ell(x)}{2} \|\theta_1 - \theta_2\|,
\end{align*}
and analogously,
\begin{align*}
1 - \sqrt{\frac{p_{\theta_1}}{p_{\theta_2} }} \leq \frac{\ell(x)}{2} \|\theta_1 - \theta_2\|.
\end{align*}
Let $A$ denote the set over which $p_{\theta_1} \geq p_{\theta_2}$, then squaring and integrating we obtain that for some sufficiently large constant $C > 0$,
\begin{align*}
H^2(p_{\theta_1},p_{\theta_2}) &\leq \|\theta_1 - \theta_2\|^2 
\left[ \int_{A} \ell^2(x) p_{\theta_1} dx + \int_{A^c} \ell^2(x) p_{\theta_2} dx \right]\\
& \leq C \|\theta_1 - \theta_2\|^2,
\end{align*}
where the final inequality uses the condition that $\sup_{\theta} P_{\theta} \ell^2 < \infty$. 

\vspace{.2cm}

\noindent {\bf Proof of Claim 2: } Fix a small $\varepsilon > 0$, by compactness of the set $\{\theta: \|\theta - \theta^*\| \geq \varepsilon\}$, and
identifiability of $\theta^*$, we obtain that, 
\begin{align*}
\inf_{\theta: \|\theta - \theta^*\| \geq \varepsilon} H(p_{\theta^*}, p_{\theta}) > 0.
\end{align*}
This in turn implies that if $H(p_{\theta^*}, p_{\theta})$ is sufficiently small, then it must be the case that $\|\theta - \theta^*\| \leq \varepsilon$. 
Locally, we can use DQM. Formally, we know that,
\begin{align*}
\int (\sqrt{p_{\theta^* + h}} - \sqrt{p_{\theta^*}} - \frac{1}{2} h^T s(\theta^*) \sqrt{p_{\theta^*}})^2 = o(h^2). 
\end{align*} Let us denote, 
\begin{align*}
\delta(x) = \sqrt{p_{\theta^* + h}} - \sqrt{p_{\theta^*}} - \frac{1}{2} h^T s(\theta^*) \sqrt{p_{\theta^*}}.
\end{align*}
Then we have that,
\begin{align*}
H^2(p_{\theta^*+h}, p_{\theta^*}) &= \int \delta^2(x) + \int \delta(x) h^T s(\theta^*) \sqrt{p_{\theta^*}} \\&~~~ + h^T I(\theta^*) h/4 \\
&\geq h^T I(\theta^*) h/4 + \int \delta^2(x) \\
&~~~- \sqrt{\int \delta^2(x)} \sqrt{ h^T I(\theta^*) h}. 
\end{align*}
Now, using the fact that $I(\theta^*)$ is non-degenerate, and that $\int \delta^2(x) = o(\|h\|^2)$ 
we obtain that for $\|h\|$ smaller than some constant,
\begin{align*}
H^2(p_{\theta^* + h}, p_{\theta^*}) \geq O(\|h\|^2).
\end{align*}
This in turn shows us that within a small ball around $\theta^*$, if $H(p_{\theta},p_{\theta^*}) \leq c_1$ (for a sufficiently small value $c_1 > 0$) 
then $H(p_{\theta}, p_{\theta^*}) \geq c_2 \|\theta - \theta^*\|$. 

\subsubsection{Proof of Lemma~\ref{lem:local_entropy}}
\begin{enumerate}
\item First let us consider $s \leq c_0$ for some small universal constant $c_0$. Using Lemma~\ref{lem:hell_main}, we have that for distributions such that $H^2(p_{\theta^*}, p_{\theta}) \leq s^2$,
it must be the case that $\|\theta - \theta^*\|^2 \leq C s^2$, for some sufficiently large universal constant $C > 0$.

As a consequence of the calculation in Example 19.7 of \cite{van2000asymptotic} we obtain that in this case for some universal constant $K > 0$,
\begin{align*}
N(u/c_3,\mathcal{F} \cap \left\{ H^2(p_\theta, p_{\theta^*}) \leq s^2 \right\} ) \leq K \left( \frac{s}{u} \right)^d.
\end{align*}

Now, integrating this we obtain that it is sufficient if:
\begin{align*}
\int_{0}^{Cs} \sqrt{ d \log (s/u)} du \leq c \sqrt{n} s^2.
\end{align*}
This is true provided that $s \geq C/\sqrt{n}$ for a sufficiently large constant $C > 0$.
\item When $s \geq c_0$, we no longer have a lower bound on the Hellinger distance in terms of the parameter distance. However, in this case a crude bound suffices. We simply bound the global metric entropy using the 
fact that $\Theta$ is compact, and once again following the calculation in Example 19.7 of \cite{van2000asymptotic}, we obtain that:
\begin{align*}
N(u/c_3,\mathcal{F} ) \leq K \left( \frac{\text{diam}(\Theta)}{u} \right)^d,
\end{align*}
and integrating this we see that it is sufficient if:
\begin{align*}
s \sqrt{d \log(1/s)} \leq c \sqrt{n} s^2,
\end{align*}
which is of course true in this regime since $s \geq c_0 > 0$.
\end{enumerate}

\subsection{Proof of Claim~\eqref{eqn:claim2}}
We use the fact that conditioned on our event we know that, $\|\widehat{\theta} - \theta^*\| \leq M/\sqrt{n}$. From Lemma 19.31 of \cite{van2000asymptotic} we obtain that,
\begin{align*}
\left|  \log \frac{\mathcal{L}_0(\widehat{\theta})}{\mathcal{L}_0(\theta^*)} - h^T G_n s_{\theta^*} + \frac{1}{2} h^T I_{\theta^*} h   \right| = o_p(1),
\end{align*}
where $G_n = \sqrt{n}(P_n - P_{\theta^*})$, $P_n$ denotes the empirical distribution on $\mathcal{D}_0$ (the first half of our samples), 
and $h = \sqrt{n} (\widehat{\theta} - \theta^*)$.
This gives us the bound:
\begin{align*}
\log \frac{\mathcal{L}_0(\widehat{\theta})}{\mathcal{L}_0(\theta^*)} =  - h^T G_n s_{\theta^*} + \frac{1}{2} h^T I_{\theta^*} h + o_p(1),
\end{align*}
where $\|h\| = O(1)$. It thus suffices to argue that, $- h^T G_n s_{\theta^*} + \frac{1}{2} h^T I_{\theta^*} h = O_p(1)$. The second term is clearly $O(1)$. 
For the first term, we apply Chebyshev's inequality. 
It is sufficient to bound the variance of $G_n s_{\theta}^*$ which is simply $I_{\theta^*}$, to obtain that, $|- h^T G_n s_{\theta^*}| = O_p(h^T I_{\theta^*} h) = O_p(1)$ 
as desired.

\section{Proof of Theorem~\ref{thm:nonpar}}
We once again build directly on Theorem~2 of \cite{wong1995}. Let $d_n$ denote the Hellinger diameter of the split LRT set. 
Once again we fix $\kappa > 0$. We want to show that, for some finite $M > 0$, 
\begin{align*}
P(d_n \geq M(\epsilon_n + \eta_n + \sqrt{\log(1/\alpha)/n})) \leq \kappa,
\end{align*}
for all $n$ large enough. 
Let us condition on the event that $\chi^2(p^*, \widehat{p}_1) \leq C_1 \eta_n^2$ throughout the proof, which holds with probability 
at least $1 - \kappa/2$ for sufficiently large $C_1 > 0$ by our assumptions. 

Observe that in our previous decomposition~\eqref{eqn:main}, we need to show that for all $p$ sufficiently far from $p^*$:
\begin{align}
\label{eqn:tt}
\log \frac{\mathcal{L}_0(p)}{\mathcal{L}_0(p^*)} + \log \frac{\mathcal{L}_0(p^*)}{\mathcal{L}_0(\widehat{p}_1)} \leq  \log(1/\alpha).
\end{align}
Theorem~2 of \cite{wong1995} guarantees us that uniformly over all $p$ such that 
$H(p,p^*) \geq \epsilon_n$, for constants $c_1, c_2 > 0$, 
\begin{align*}
\log \frac{\mathcal{L}_0(p)}{\mathcal{L}_0(p^*)} \leq - c_1 n H^2(p,p^*),
\end{align*} 
with probability at least $1 - c_2 \exp(- n H^2(p,p^*)).$ For the second term we observe that for any $C > 0$, 
\begin{align*}
P\left( \frac{\mathcal{L}_0(p^*)}{\mathcal{L}_0(\widehat{p}_1)} \geq \exp(C n s^2)\right) &\leq \mathbb{E} \left[ \frac{\mathcal{L}_0(p^*)}{\mathcal{L}_0(\widehat{p}_1)}\right] \exp(- Cn s^2) \\
&= (1 + \chi^2(p^*, \widehat{p}_1))^n \exp(- Cn s^2) \\
& \leq \exp(n C_1 \eta_n^2) \exp( - Cn s^2).  
\end{align*}
Putting these two results together, we see that with probability at least $1 - \kappa/2$ for any distribution $p \in \mathcal{P}$ if $H(p,p^*) \geq M(\epsilon_n + \eta_n + \sqrt{\log(1/\alpha)/n})$ for a sufficiently large constant $M > 0$,
then~\eqref{eqn:tt} is satisfied and $p$ will not belong to our confidence set. Thus, we obtain the desired diameter bound.

\end{document}